\documentstyle[12pt]{article}
\begin{document}

\newcommand\comp{\circ}
\newtheorem{theorem}{Theorem}[section]
\newtheorem{example}[theorem]{Example}
\newtheorem{conjecture}[theorem]{Conjecture}
\newtheorem{examples}[theorem]{Examples}
\newtheorem{proposition}{Proposition}[section]
\newtheorem{fact}[theorem]{Fact}
\newtheorem{problem}[theorem]{Open problem}

\newcommand{\I}{\mbox{{\bf I}}}
\newcommand{\A}{\mbox{{\bf Ax \ }}}
\newcommand{\qed}{\hspace*{\fill}\rule{1 ex}{1.5 ex}\\}

\newtheorem{corollary}[theorem]{Corollary}
\newtheorem{definition}[theorem]{Definition}
\newtheorem{remark}[theorem]{Remark}
\newtheorem{lemma}[theorem]{Lemma}
\newtheorem{claim}{Claim}[theorem]
\newtheorem{construction}[theorem]{Construction}
\newenvironment{proof}{\begin{trivlist}\item[]{\bf
Proof.}}{\qed\end{trivlist}}

\newcommand{\comment}[1]{\typeout{swallowing comment}}

\hyphenation{equi-v-al-ent dimen-s-ion mon-ad-ic}

\title{On topological properties of ultraproducts of finite sets}
\author{G\'abor S\'agi\thanks{Supported by Hungarian
National Foundation for Scientific Research grant D042177.},
Saharon Shelah\footnote{The second author would like to thank the
Israel Science Foundation for partial support of this research
(Grant no. 242/03). Publication 841.}}


\maketitle
\begin{abstract}
Motivated by the model theory of higher order logics, in
\cite{cl1a} a certain kind of topological spaces had been introduced on
ultraproducts. These spaces are called ultratopologies. Ultratopologies
provide a natural extra topological structure for ultraproducts and using
this extra structure in \cite{cl1a} some preservation and characterization
theorems had been obtained for higher order logics. \\
\indent
The purely topological properties of ultratopologies seem interesting on their
own right. We started to study these properties in \cite{ut}, where some
questions remained open. Here we present the solutions of two such problems.
More concretely we show that \\
\indent $(1)$ there are sequences of finite sets of pairwise
different cardinalities such
that in their certain ultraproducts there are homeomorphic ultratopologies and \\
\indent $(2)$ if $A$ is an infinite ultraproduct of finite sets
then every ultratopology on $A$ contains a dense subset $D$ such
that $|D| < |A|$.
\\
\\
{\em AMS Classification: } 03C20, 54A25, 54A99.  \\
{\em Keywords:} ultraproduct, ultratopology, dense set.
\end{abstract}

\section{Introduction}
\label{intro}

In first order model theory the ultraproduct construction can be
applied rather often. this is because ultraproducts preserve the
validity of first order formulas. It is also natural to ask, what
connections can be proved between certain higher
order formulas and ultraproducts of models of them. \\
\indent
In \cite{cl1a} we answer related questions in terms of topological
spaces which can be naturally associated to ultraproducts. These spaces are
called ultratopologies and their definition
can be found in \cite{cl1a} and also at the beginning of
\cite{ut}. \\
\indent Although ultratopologies were introduced from
logical (model theoretical) reasons, these spaces can be
interesting on their own right. In \cite{ut} a systematic
investigation about these topological properties has been started.
However, in \cite{ut} some problems remained open. In the present note we
are dealing with two such problems. \\
\indent
In Section \ref{hb} we give a positive answer for Problem 5.2 of
\cite{ut}: there are sequences of finite sets of pairwise
different cardinalities such that their certain ultraproducts are
still homeomorphic with respect to some carefully chosen
ultratopologies. In fact, in Theorem \ref{discr} below we show
that if $U$ is a good ultrafilter and $A$ is any infinite
ultraproduct of finite sets modulo $U$ then there is an
ultratopology on $A$ in which the family of closed sets consists
just the finite sebsets of $A$ and the whole $A$. From this the
affirmative answer for Problem 5.2 of \cite{ut} can be
immediately deduced. In Section \ref{ds} we investigate the
possible cardinalities of dense sets in ultratopologies, again on
ultraproducts of finite sets. In Corollary \ref{dc} we show that
if ${\cal C}$ is any ultratopology on an infinite ultraproduct
$A$ then the density of ${\cal C}$ is smaller then $|A|$, that
is, one can always find a
dense set whose cardinality is less than $|A|$. \\
\indent
Throughout we use the following conventions. $I$ is a set
and for every $i \in I$ $A_{i}$ is a set. Moreover $ A = \Pi_{i
\in I} A_{i} / U$ denotes the ultraproduct of $A_{i}$'s modulo an
ultrafilter $U$. \\
\indent Every ordinal is the set of smaller ordinals and natural
numbers are identified with finite ordinals. Throughout $\omega$
denotes the smallest infinite ordinal
and $cf$ denotes the cofinality operation. \\
\indent
In order to simplify notation, sometimes we will identify
${}^{k}(\Pi_{i \in I}A_{i})$ with $\Pi_{i \in I}{}^{k}A_{i}$ by
the natural way, that is, $k$--tuples of sequences
are identified with single sequences whose terms are $k$--tuples. \\
\indent
If
$X$ is a topological space and $A \subseteq X$ then $cl(A)$
denotes the closure of $A$. Suppose $k \in \omega$, $\langle
A_{i}: \ i \in I \rangle$ is a sequence of sets, $U$ is an
ultrafilter on $I$ and $R_{i} \subseteq {}^{k}A$ is a given
relation for every $i \in I$. Then the {\em ultraproduct} relation
$\Pi_{i \in I} R_{i} /U$ is defined as follows. \\
\\
\centerline{$ \Pi_{i \in I} R_{i} /U = \{ s/U \in ({}^{k} \Pi_{i \in I} A_{i}/U): \  \{ i \in I: \ s(i) \in R_{i}   \} \in U  \}$.} \\
\\
As we mentioned, we assume that the reader is familiar with the
notions of ``choice function'', ``ultratopology'', ``a point is
close to a relation'', ``$T(a,R)$'', etc. These notions were
introduced in \cite{cl1a} and a short (but fairly complete)
survey can be found at the beginning of \cite{ut}.


\section{Homeomorphisms between different ultraproducts}
\label{hb}

In \cite{ut} Problem 5.2 asks whether is it possible to choose ultrafilters
$U_{1}$, $U_{2}$ and sequences of natural numbers
$s = \langle n_{i}, \ i \in I \rangle$ and
$ z = \langle m_{j}, \ j \in J \rangle$ so that \\
\indent $(i)$ $n_{i} \not= m_{j}$ for all $i \in I, j \in J$ and \\
\indent $(ii)$ for every $k \in \omega$ there are $k$--dimensional ultratopologies
${\cal C}_{k}$ in $\Pi_{i \in I} n_{i} / U_{1}$ and ${\cal D}_{k}$ in
$\Pi_{j \in J} m_{j} / U_{2}$ such that ${\cal C}_{k}$ and ${\cal D}_{k}$ are
homeomorphic?

We will give an affirmative answer. In fact, we prove the
following theorem from which the above question can be easily
answered.

\begin{theorem}
\label{discr} Suppose $\langle n_{i}, \ i \in I \rangle$ is an
infinite sequence of natural numbers and $U$ is a good
ultrafilter on $I$ such that $A = \Pi_{i \in I} n_{i} / U$ is
infinite. Then for every $k \in \omega$ there is a
$k$--dimensional choice function on $A$ such that the family of
closed sets in the induced ultratopology consists of the finite
subsets of ${}^{k}A$ and ${}^{k}A$.
\end{theorem}

\begin{proof}
Let $c$ be any $k$--dimensional choice function on ${}^{k}A$. By modifying $c$,
we will construct another choiche function $\hat{\ }$ which induces the required
ultratopology.
Let $E$ be the set of all triples $ \langle s, i, m \rangle$ where
$s \in \Pi_{l \in I} {\cal P}({}^{k}n_{l})$, $i \in I$ such that
$\Pi_{l \in I} s(l) / U$ is infinite and $m \in {}^{k}n_{i} -
s(i)$. It is easy to see that $|\Pi_{i \in I} {\cal
P}({}^{k}n_{i})| \leq |{}^{I} \omega| = 2^{|I|}$. Therefore $|E|
\leq 2^{|I|} \times |I| \times \aleph_{0} = 2^{|I|}$. Let $\{
\langle s_{\alpha}, i_{\alpha}, m_{\alpha} \rangle: \ \alpha <
2^{|I|} \}$
be an enumeration of $E$. \\
\indent By transfinite recursion we construct an injective
function $f: E \rightarrow {}^{k}A$ such that for every $\langle
s, i, m \rangle \in E$ one has $f(\langle s, i, m \rangle) \in
\Pi_{l \in I} s(l) / U$. Suppose $f_{\beta}$ has already been
defined on $ \{ \langle s_{\alpha}, i_{\alpha}, m_{\alpha}
\rangle: \ \alpha < \beta \}$
for all $\beta < \gamma \leq 2^{|I|}$ such that \\
\indent $\beta_{1} < \beta_{2} < \gamma \Rightarrow f_{\beta_{1}} \subseteq f_{\beta_{2}}$ and \\
\indent $|f_{\beta}| \leq |\beta|$. \\
If $\gamma$ is a limit ordinal, then let $f_{\gamma} =
\cup_{\beta < \gamma} f_{\beta}$. Now suppose $\gamma$ is a
successor ordinal, say $\gamma = \delta + 1$. Since $U$ is a good
ultrafilter, by Theorem VI, 2.13 of \cite{clasth} it follows,
that the cardinality of $B = \Pi_{l \in I} s(l) / U $ is
$2^{|I|}$. Therefore there is an element $b \in B$ which is not
in the range of $f_{\delta}$. Let $f_{\gamma}$ be $f_{\delta}
\cup \{ \langle \langle s_{\delta}, i_{\delta}, m_{\delta}
\rangle,  b \rangle \}$.
Clearly, $f = f_{2^{|I|}}$ is the required function. \\
\indent
Now we construct a $k$--dimensional choice function $\hat{\ }$ as follows.
If $a \not\in rng(f)$ then let $\hat{a} = c(a)$. Otherwise there is a
unique $ e = \langle s, i, m \rangle \in E$ such that $f(e) = a$.
Let
\[  \hat{a}(j) =  \left\{ \begin{array}{ll}
                                 c(a)(j) & \mbox{if $i \not= j$, } \\
                                 m  &   \mbox{otherwise.}
                            \end{array}
                   \right. \]
\\
In this way we really defined a $k$--dimensional choice function
on $A$. We claim that the closed sets of the induced
ultratopology are exactly the finite subsets of ${}^{k}A$
and ${}^{k}A$. \\
\indent
By Theorem 2.5 of \cite{ut} every ultratoplogy is $T_{1}$ therefore
every finite subset of ${}^{k}A$ is closed. Let $F$ be an infinite closed
subset of ${}^{k}A$ and suppose, seeking a contradiction, that there is an
element \\
\\
\centerline{$(*)$ \indent $b \in {}^{k}A - F$.} \\
\\
By Corollary 2.2 of \cite{ut} $F$ is a decomposable relation, say
$F = \Pi_{l \in I} s(l) / U$. Therefore, there is a $J \in U$ such
that for every $i \in J$ one has $\hat{b}(i) \not\in s(i)$.
Hence, for every $i \in J$ $e_{i} = \langle s, i, \hat{b}(i)
\rangle \in E$.
By construction, for every $i \in J$ one has $f(e_{i}) \in F$ and
$f(e_{i})\hat{\ }(i) = \hat{b}(i)$.
This means that \\
\\
\centerline{ $\{ i \in I: \exists a \in F: \hat{a}(i) = \hat{b}(i) \} \supseteq J \in U$.} \\
\\
That is, $T(F,b) \in U$ (where $T$ is understood according to the
new choice function $\hat{\ }$). Since we assumed that $F$ is
closed, this implies $b \in F$ which contradicts to $(*)$.
\end{proof}

\begin{corollary}
There are ultrafilters $U_{1}$, $U_{2}$ (respectively, over $I$
and $J$) and sequences of natural numbers $s = \langle n_{i}, \ i
\in I \rangle$ and
$ z = \langle m_{j}, \ j \in J \rangle$ so that \\
\indent $(i)$ $n_{i} \not= m_{j}$ for all $i \in I, j \in J$ and \\
\indent $(ii)$ for every $k \in \omega$ there are $k$--dimensional ultratopologies
${\cal C}_{k}$ in $\Pi_{i \in I} n_{i} / U_{1}$ and ${\cal D}_{k}$ in
$\Pi_{j \in J} m_{j} / U_{2}$ such that ${\cal C}_{k}$ and ${\cal D}_{k}$ are
homeomorphic.
\end{corollary}

\begin{proof}
Let $U_{1}, U_{2}$ be good ultrafilters and let $s$ and $z$ be arbitrary
sequences of natural numbers satisfying the requirements
of the corollary such that $|I| = |J|$ and both $A = \Pi_{i \in I} s_{i} / U_{1}$ and
$B = \Pi_{j \in J} z_{j} / U_{2}$ are infinite. Let $ k \in \omega$ be arbitrary.
By Theorem \ref{discr} there are ultratopologies ${\cal C}_{k}$, ${\cal D}_{k}$,
respectively on $A$ and $B$ such that \\
\indent the closed sets of ${\cal C}_{k}$ are exactly the finite subsets of
${}^{k}A$ and ${}^{k}A$ and \\
\indent the closed sets of ${\cal D}_{k}$ are exactly the finite subsets of
${}^{k}B$ and ${}^{k}B$. \\
By Theorem VI, 2.13 of \cite{clasth} $|A|=|B| = 2^{|I|}$. Let $f: A \rightarrow B$
be any bijection. Then $f_{k}: {}^{k}A \rightarrow {}^{k}B$,
$f(\langle a_{0},...,a_{k-1} \rangle = \langle f(a_{0}),...,f(a_{k-1}) \rangle$
is clearly a bijection from ${}^{k}A$ onto ${}^{k}B$ mapping finite subsets
of ${}^{k}A$ to finite subsets of ${}^{k}B$. Thus, $f$ is the required
homeomorphism.
\end{proof}

\section{cardinalities of dense sets}
\label{ds}

Problem 5.3 (A) of \cite{ut} asks whether is it possible to
choose a sequence of finite sets $s$ and an ultrafilter $U$ so
that there is an utratopology on $A = \Pi_{i \in I} n_{i} / U$ in
which every dense set has cardinality $|A|$. In this section we
will show that this is impossible if $A$ is infinite. We start by
a simple observation: every $k$-dimensional ultratopology is
homeomorphic with an appropriate $1$--dimensional ultratopology.

\begin{theorem}
\label{dims} Suppose ${\cal C}_{k}$ is a $k$--dimensional
ultratopology on $A$. Then there is
a $1$--dimensional ultratopology ${\cal D}$
which is homeomorphic to ${\cal C}$.
\end{theorem}

\begin{proof}
The idea is to identify $k$--tuples of sequences by sequences of
$k$--tuples. By a slight abuse of notation, we will use this
identification freely. Let $A = \Pi_{i \in I} A_{i} / U$ (here
the $A_{i}$'s are not necessarily finite) and suppose $\hat{\ }$
is a $k$--dimensional choice function inducing ${\cal C}_{k}$.
Let ${\cal B} = \Pi_{i \in I} {}^{k}A_{i} / U$. We define a
$1$--dimensional choice function $c$ in $B$ as follows. If
$s=\langle s_{i}: i \in I \rangle / U \in B$ then for each $j \in
k$ let $s^{j} = \langle s_{i}(j): \ i \in I \rangle / U$. Define
$c(s) = \langle s^{0},...,s^{k-1} \rangle\hat{\ }$ and $\varphi:
{}^{k}A \rightarrow B$, $\varphi (\langle s^{0}/ U,...,s^{k-1}/U
\rangle) = \langle \langle s^{0}(i),...,s^{k-1}(i) \rangle: i \in
I \rangle / U$. Then clearly, $c$ is a $1$--dimensional choice
function which induces an ultratopology ${\cal D}$ on $B$. Then
for any $a \in {}^{k}A$ and $i \in I$ one has $\hat{a}(i) =
c(\varphi(a))(i)$. Now it is straightforward to check that
$\varphi$ is a homeomorphism between ${\cal C}$ and ${\cal D}$.
\end{proof}

Let ${\cal C}$ be an ultratopology on an ultraproduct $A = \Pi_{i \in I} n_{i} / U$
of finite sets. Suppose ${\cal C}$ can be induced by a choice function $\hat{\ }$.
Let \\
\\
\centerline{$G = \{ \langle i,m \rangle: \ i \in I, m \in n_{i}$ and $( \exists a \in A)(\hat{a}(i) = m) \}$} \\
\\
and for each $\langle i,m \rangle \in G$ let $a_{i,m} \in A$ be
such that $\hat{a}_{i,m}(i)=m$. Clearly, if $I$ is infinite, then
$|G| \leq |I|$. We claim that there is a dense subset $R$ of $A$
such that $|R| \leq |G|$. In fact, $R$ can be chosen to be $R =
\{ a_{i,m}: \ \langle i,m \rangle \in G \}$. To see this, suppose
$a \in A$. Then for every $i \in I$ one has $\langle i,
\hat{a}(i) \rangle \in G$ and therefore $T(R,a) = I \in U$. Hence
$cl(R) = A$, as desired. \\
\indent
Now we are able to provide a negative answer for Problem 5.3 (A) of \cite{ut}.

\begin{theorem}
\label{maint} Suppose ${\cal C}$ is a $1$--dimensional
ultratopology on an infinite ultraproduct $A = \Pi_{i \in I}
n_{i} / U$ where each $n_{i}$ is a finite set. Then $d({\cal C})
< |A|$.
\end{theorem}

\begin{proof}
Suppose, seeking a contradiction, that ${\cal C}$ is an
ultratopology on $A$ such that the cardinality of every dense set
in ${\cal C}$ is equal with $|A| = \kappa$. Using the notation
just introduced in the remark before the theorem,
$R$ is a dense subset of $A$ and therefore $|R| = |A|$.
Let $<^{A}$ be a well ordering of $A$ (having order type $\kappa$).
By transfinite recursion we define a sequence $\langle b_{i} \in
A: i < \kappa \rangle$ as follows. Assume $j < \kappa$ and
$b_{l}$ has already been defined for every $l < j$. Let $W_{j} =
\{ b_{l}: l < j \}$ and let $V_{j} = \{ a \in A: T(W_{j},a) \in U
\}$.
Since $j < \kappa$, $V_{j} \subseteq cl(W_{j}) \not= A$.
If $j$ is an odd ordinal,
then let $b_{j}$ be
the $<^{A}$--first element of $R- W_{j}$. Otherwise let $b_{j}$ be the
$<^{A}$--first element in $A-V_{j}$. Clearly, the following conditions are
satisfied: \\
\indent $(i)$ for every $\langle i,m \rangle \in G$ there is a $j
< \kappa$ such that $\hat{b}_{j}(i) = m$,
in fact, $R \subseteq \{b_{l}: \ l < \kappa \} = W_{\kappa}$. \\
\indent $(ii)$ for every $a \in A$ there is a $j < \kappa$ such that $a \in V_{j}$
(the smallest such $j$ will be denoted by $j_{a}$), \\
\indent $(iii)$ for every $j < \kappa$ there is an ordinal $j < s(j)$ such that
$s(j) < \kappa$ and $b_{s(j)} \not\in cl(W_{j})$. (This is true because otherwise
by $(i)$ one would have $cl(W_{j}) \supseteq cl(R) = A$ which is impossible
since $|W_{j}| < \kappa$.) \\
Now let $H = \{ i \in I: n_{i} = \{ \hat{b}_{j}(i): j < \kappa \}
\}$. We show that $H \in U$. \\
\indent Again, seeking a contradiction, assume $I-H \in U$. For
every $i \in I-H$ let $e_{i} \in n_{i} - \{ \hat{b}_{j}(i): j <
\kappa \}$ be arbitrary, let $e = \langle e_{i}: i \in I-H
\rangle / U$ and let $O = \{i \in I: \hat{e}(i) = e_{i} \}$.
Clearly, $O \in U$. In addition, if $i \in O \cap (I-H)$ then by
$(i)$ there is a $j < \kappa$ such that $\hat{b}_{j}(i) =
\hat{e}(i) = e_{i}$
contradicting to the selection of $e_{i}$'s. \\
\indent For every $i \in H$ we introduce a binary relation
$\prec^{i}$ as follows. If $n,m \in n_{i}$ then $n \prec^{i} m$
means that there is a $j \in \kappa$ such that $\hat{b}_{j}(i) =
n$ but for every $l \leq j$ $\hat{b}_{l}(i) \not= m$. Clearly,
for each $i \in H$ the relation $\prec^{i}$ is irreflexive,
transitive, and trichotome. Since $n_{i}$ is finite, for every $i
\in H$ there is an $\prec^{i}$--maximal element $y_{i} \in
n_{i}$. Let $y = \langle y_{i}: \ i \in H \rangle / U$ and let $K
= \{ i \in I: \hat{y}(i) = y_{i} \}$. Now by $(ii)$ and $(iii)$
$y \in V_{j_{y}}$ and $b_{s(j_{y})} \not\in cl(W_{j_{y}}) =
cl(V_{j_{y}})$.
Thus, for every $i \in K \cap T(W_{j_{y}},y)$ one has \\
\\
\centerline{$ (*) \indent \hat{y}(i) = y_{i} \in \{ \hat{b}_{l}(i): l < j_{y} \}.$} \\
\\
$cl(V_{j_{y}})$ is closed, therefore there is $L \in U$ such that
for all $i \in L$ one has $\hat{b}_{s(j_{y})}(i) \not \in \{
\hat{b}_{l}(i): l < j_{y} \}$ and thus $\hat{b}_{l}(i) \prec^{i}
\hat{b}_{s(j_{y})}(i)$ for every $i \in L \cap H$ and for every
$l < j_{y}$. Particularly, $(*)$ implies that if $i \in L \cap H
\cap K \cap T(W_{j_{y}},y)$ then $\hat{y}(i) = y_{i} \prec^{i}
\hat{b}_{s(j_{y})}(i)$ which is impossible since by construction,
$y_{i}$ is the $\prec^{i}$--maximal element in $n_{i}$. This
contradiction completes the proof.
\end{proof}

Using Theorem \ref{dims} the above results can be generalized to higher
dimensional ultratopologies as well.

\begin{corollary}
\label{dc}
Let $k \in \omega$ be arbitrary and suppose ${\cal C}$ is a $k$--dimensional
ultratopology on an infinite ultraproduct
${\cal A} = \Pi_{i \in I} n_{i} / U$ where each $n_{i}$ is
a finite set. Then $d({\cal C}) < |A|$.
\end{corollary}

\begin{proof}
Assume, seeking a contradiction, that ${\cal C}$ is a $k$--dimensional
ultratopology on $A$ whose density is $|A|$. Then, by Theorem \ref{dims}
there is a $1$--dimensional ultratopology with the above property, contradicting
to Theorem \ref{maint}
\end{proof}

\bigbreak
\leftline{Alfr\'ed R\'enyi Institute of Mathematics}
\leftline{Hungarian Academy of Sciences}
\leftline{Budapest Pf. 127}
\leftline{H-1364 Hungary}
\leftline{sagi@renyi.hu}

\ \\

\bigbreak
\leftline{Department of Mathematics}
\leftline{Hebrew University}
\leftline{91904 Jerusalem, Israel}
\leftline{shelah@math.huji.ac.il}


\begin{thebibliography}{99}

\bibitem {chk} {\sc C.C. Chang, H.J. Keisler,}
{\it Model Theory, \/} North--Holland, Amsterdam (1973).
\bibitem {cl1a} {\sc G. S\'agi,}
{\it Ultraproducts and higher order formulas, \/}
Math. Logic Quarterly, Vol. 48, No. 2, pp. 261--275, (2002).
\bibitem {ut} {\sc J. Gerlits, G. S\'agi,}
{\it Ultratopologies, \/} Accepted for publocation, Math. Logic
Quarterly, (2004).
\bibitem {clasth} {\sc S. Shelah,}
{\it Classification theory, \/} North--Holland, Amsterdam (1990).
\end{thebibliography}
\end{document}